\documentclass[12pt]{article}
\usepackage{mathrsfs}
\usepackage{amsmath,amsfonts,amssymb,rotating,amsthm}
\usepackage{bm}
\usepackage{hyperref}
\usepackage{array}
\usepackage{breqn}
\usepackage{enumitem}
\allowdisplaybreaks[4]
\def\qed{\nopagebreak\hfill{\rule{4pt}{7pt}}}
\def\proof{\noindent {\it{Proof.} \hskip 2pt}}
\newtheorem{theo}{Theorem}[section]
\newtheorem{deff}[theo]{Definition}
\newtheorem{lemm}[theo]{Lemma}

\newtheorem{conj}[theo]{Conjecture}

\theoremstyle{remark}
\newtheorem{remark}[theo]{Remark}

\numberwithin{equation}{section}

\newdimen\Squaresize \Squaresize=11pt
\newdimen\Thickness \Thickness=0.7pt
\def\Square#1{\hbox{\vrule width \Thickness
   \vbox to \Squaresize{\hrule height \Thickness\vss
    \hbox to \Squaresize{\hss#1\hss}
   \vss\hrule height\Thickness}
\unskip\vrule width \Thickness} \kern-\Thickness}

\def\Vsquare#1{\vbox{\Square{$#1$}}\kern-\Thickness}

\def\moins{\raise 1pt\hbox{{$\scriptstyle -$}}}

\begin{document}

\begin{center}
{\Large \bf Finding Product and Sum Patterns in non-commutative settings}
\end{center}

\begin{center}
T.Y. Tao$^{1}$ and Neil N.Y. Yang$^{2}$\\[8pt]
$^{1,2}$School of Mathematical Sciences\\
Fudan University\\
Shanghai 200433, P. R. China\\[6pt]
Email: {\tt $^{1}$tytao20@fudan.edu.cn, $^{2}$23110180046@m.fudan.edu.cn \tt}
\end{center}

\vspace{0.3cm} \noindent{\bf Abstract.} Hindman conjectured that any finite partition of $\mathbb{N}$ has a monochromatic $\{x,y,x+y,xy\}$. Recently, Bowen proved the result for all 2-partition. In this paper, we extend Bowen's result to any semiring $(S,+,\cdot)$ such that $Ss$ is piecewise syndetic for all $s\in S$. As a method, we gave a combinatorial proof for a piecewise syndetic version of Bergerson and Glasscock's {\bf IP}$_r^*$ Szemer\'edi Theorem, and discussed the case when the operation is not commutative.
\\

\\

\section{Introduction}

Partition regularity is an important field in additive combinatorics, which is old but still full of vitality. It dates back to Roth's Theorem and van der Waerden's Theorem, which states that one can arithmetic progressions in any finite partition of $\mathbb{N}$. These results are quickly extended to arbitrary semigroups and so are a lot of similar results.

One of the most famous conjectures in additive combinatorics, which is conjectured by Hindman, asks whether for any finite partition of $\mathbb{N}$, the set of positive integers, there exists a monochromatic $\{x,y,x+y,xy\}$ (i.e. one can find a subset in the partition, and elements x and y, such that all of $\{x,y,x+y,xy\}$ are in the subset). He proved in \cite{hindman} that any 2-coloring of $\{2,\dots,990\}$ contains a monochromatic set of such form. The conjecture is still open and many works have been done recently. Moreira \cite{moreira} proved that one can find monochromatic $\{x,x+y,xy\}$ in any finite partition of $\mathbb{N}$ (and also similar patterns in LIDs). Bowen solved a special case when the partition has size 2. The following theorem is proved by Bowen \cite{bowen}.

\begin{theo}[Bowen]
Given any $2$-coloring of $\mathbb{N}$ and $n\in\mathbb{N}$, there are arbitrarily
large and distinct $x_1, ..., x_n\in\mathbb{N}$ such that
\[
\{x_i,\prod_{j\leq i}x_j,\sum_{i=1}^nx_j:i\leq n\}
\]
is monochromatic. Moreover, there are also monochromatic sets of the form $\{x,y,xy,x+ny\}$.
\end{theo}

More recently, Bowen and Sabok \cite{bs} and Alweiss \cite{alweiss} independently proved that one can find monochromatic $\{x,y,x+y,xy\}$ in $\mathbb{Q}$, which turns out to be very different from $\mathbb{N}$ since the multiplication in $\mathbb{Q}$ has its inverse. They used different methods and obtained different generalizations of Hindman's Conjecture on $\mathbb{Q}$.

One of the aims of our paper is to extend Bowen's result on $\mathbb{N}$ to semirings that are additively large in some way, as well as giving some ideas and useful tools on how do do similar problems without the commutative law. Recall that a semiring $(S,+,\cdot)$ is a set $S$ and two binary operations, with the commutative law in $+$ and the distributive law. We state one of our theorem as follows.

\begin{theo}
For any positive integer $k$ and any $2$-coloring of a semiring $S=(S,+,\cdot)$ such that $Ss$ is piecewise syndetic in $(S,+)$ for all $s\in S$, there is either a monochromatic $\{x,y,x+y,x^iy:i\leq k\}$, or a monochromatic $\{x,y,xy,x+iy:i\leq k\}$.
\end{theo}

Recall that $A\subset (S,\circ)$ is thick if $\forall k\in \mathbb{N}$ and $s_1,\dots,s_k\in S,\exists x\in S$ such that $\{x\circ s_1,\dots,x\circ s_k\}\subseteq A$, and is (right) multiplicatively syndetic if $\exists k\in \mathbb{N}$ and $s_1,\dots,s_k\in S$ such that $S=s_1^{-1}A\cup\dots\cup s_k^{-1}\circ A$, where $x^{-1}A$ denotes the set $\{s\in S:x\circ s\in A\}$. We also discuss whether we can expect to find such a set in a sufficiently large subset of $S$.

\begin{theo}
For any semiring $S=(S,+,\cdot)$ and $A\subseteq S$, if $A$ is both multiplicatively thick and syndetic, then there is $x,y$ such that $\{x,y,xy,x+y\}\subset A$.
\end{theo}

To prove these results in semirings, we first need to generalize a variant of ${\bf IP}_r^*$ Szemer\'edi's Theorem to semirings. Recall that for a muti-subset $A$ with elements in $S$, ${\bf FS}(A)$ denotes the set of all finite sums with elements in $A$. We say that $A\subseteq S$ is ${\bf IP}_r^*$ if for any $B\subseteq S$ of $r$ elements, $A\cap {\bf FS}(B)\neq\emptyset$. 

We need the following theorem which was proved by Furstenberg and Katznelson \cite{fk}.

\begin{theo}[Furstenberg]\label{fk}
Let $(X,\mathcal{B},\mu)$ be a measure space with $\mu(X)<\infty$, let $T_1,T_2,\dots T_n$ be commuting measure preserving transformations of $X$ and let $A\in\mathcal{B}$  with $\mu(A)>0$. Then
\[ 
\liminf_{M\to\infty}\frac{1}{M}\sum_{m=1}^{M}\mu\left(\bigcap_{i=1}^kT_n^{-m}A\right)>0.
\]
\end{theo}

In particular, for every finite partition $X=\bigcup_{j=1}^{k}C_j$, there is some $j$ such that $\mu(C_j)>0$,  and hence there is some $m$ such that $\bigcap_{i=n}^kT_n^{-m}C_j\neq\emptyset$.

Bergelson and Glasscock \cite{bg} gave a short proof for the following statement. Recall that a left invariant mean on a semigroup $S$ is a positive linear functional $\lambda$ of norm $1$ which is left translation invariant, i.e. $\lambda\left(s^{-1}A\right)=\lambda(A)$ for all $s\in S$ and $A\subseteq S$. Especially, every commutative semigroup admits an invariant mean. For detailed discussions on this topic, see \cite{bg}.

\begin{theo}[Bergelson \& Glasscock]
Let $n\in\mathbb{N}$ and $\delta>0$. There exists $r\in\mathbb{N}$ and $\beta>0$ for which the following holds. For all commutative semigroups $(S,+)$ and $(R,+)$, all homomorphisms $\varphi_1,\dots, \varphi_n$ from $S$ to $R$, all invariant means $\lambda$ on $R$, and all $A\subseteq R$ with $\lambda(A)\geq\delta$, the set
\[
\left\{s\in S:\lambda\left(\bigcap_{i=1}^{n}(A-\varphi_i(s))\right)\geq\beta\right\}
\]
is ${\bf IP}_r^*$ in $S$.
\end{theo}

We shall derive the following version in piecewise syndetic version by a combinatorial method. Recall that $A\subseteq (R,\cdot)$ is piecewise syndetic if $A$ is the intersection of a thick set and a syndetic set.

\begin{theo}
Let $S=(S,\cdot)$ and $R=(R,\cdot)$ be semigroups such that $R$ admits a left invariant mean, and $A\subseteq R$ be piecewise syndetic. For homomorphisms $\varphi_1,\dots, \varphi_n$ from $S$ to $R$, define
\[
D(A;\varphi_1,\dots, \varphi_n)=\left\{d\in S:\bigcap_{i=1}^{n}(\varphi_i(d)^{-1}A)\quad \text{is piecewise syndetic}\right\},
\]
where $Aa^{-1}=\{b\in R:ba\in A\}$ Then we have
\begin{enumerate}
\item For any $n\in\mathbb{N}$ and $\varphi_1,\dots, \varphi_n$ such that $\varphi_i\varphi_j=\varphi_j\varphi_i$ for all $i,j$, there is some $t\in R$ such that $D(t^{-1}A;\varphi_1,\dots, \varphi_n)\neq\emptyset$.
\item Suppose in addition that $S$ and $R$ are commutative. Then for any $n\in\mathbb{N}$, there is $r\in\mathbb{N}$ such that for any $\varphi_1,\dots, \varphi_n$, $D(A;\varphi_1,\dots, \varphi_n)$ is ${\bf IP}_r^*$ in $S$.
\end{enumerate}
\end{theo}

In the proof of this theorem, we shall use a technique to derive a piecewise syndetic set from a single element. This technique can also be used when we extend Bowen's result to semirings, finding a piecewise syndetic set of patterns of $\{x,y,x+y,xy\}$ if the partition satisfies some property.

The rest of this paper is organized as follows. In section 2 we give some definitions and background information for largeness in semigroups, which are fundamental in the proof of our theorems. In section 3 we prove our variant of ${\bf IP}_r^*$ Szemer\'edi's Theorem, and also give a technique to find a large set of certain patterns in non-commutative setting. In section 4 we prove Bowen's Theorem in semirings, and give some generalizations in certain cases. At last we propose some conjectures relevant to our methods.

\section{Large subsets}

We are going to systematically define some kinds of largeness, some of which rely on the algebraic structure of semigroup $(\beta S,\circ)$ (here $\circ$ is any semigroup operation on $S$). Recall that $\beta S$ is the space of ultrafilters of $S$, equipped with the topology generated by the clopen sets 
\[
\widehat{A}=\{p\in\beta S:A\in p\}
\]
through all $A\subseteq S$. The ultrafilter space is only an auxiliary tool and have few appearances in our main results. One only need to bear in mind the we can extend $\circ$ on $\beta S^2$ to make $(\beta S, \circ)$ a compact right topological semigroup. If $e\in L$ where $L\subseteq\beta S$ is a minimal left ideal and $e$ is an idempotent, then we call $e$ a minimal idempotent. For a thorough treatment of the ultrafilter space and the algebra on it, see \cite{hs}.

Now we are ready to define largeness of subsets of $S$. To be concise, we denote $\mathcal{P}_f(S)$ the collection of finite muti-subsets with elements in $S$. We will see why we use mutisets when we are defining {\bf IP} system. Some of the concepts have a connection between their conbinatorical definitions and their definitions using the ultrafilter space, also see \cite{bg} for more background information.

\begin{deff}
For a semigroup $(S,\circ)$ and $A\subseteq S$, we say that
\begin{itemize}
\item $A$ is thick if $\forall F\in\mathcal{P}_f(S), \exists x\in S$ such that $F\circ x\subseteq A$. Equivalently, if there is a minimal left ideal $L\in\beta S$ such that $L\in\widehat{A}$.
\item $A$ is syndetic if $\exists F\in\mathcal{P}_f(S)$ such that $S=\bigcup_{s\in F} s^{-1}\circ A$. Equivalently, if for all minimal left ideal $L\in\beta S,L\cap\widehat{A}\neq\emptyset$.
\item $A$ is piecewise syndetic if $\exists F\in\mathcal{P}_f(S)$ such that $\bigcup_{s\in F} s^{-1}\circ A$ is thick. Equivalently, if there is a minimal left ideal $L\in\beta S$ such that $L\cap\widehat{A}\neq\emptyset$.
\item $A$ is central if there is a minimal idempotent $e\in\beta S$ such that $A\in e$.
\end{itemize}
\end{deff}
Concepts here can be defined both for left and for right, but we will only use left or right that corresponds to this definition in this paper and we will omit these attributes.

It is natural to ask whether the homomorphism image of a large set is still large. We will need the following theorem stated by Bergelson and Glasscock \cite{bg}.

The following theorem is also due to Hindman and Strauss \cite{bg}.

\begin{theo}[Hindman \& Strauss]\label{homo}
Let $(S, \cdot),(T,\cdot)$ be semigroups, $\varphi:(S,\cdot)\to(T,\cdot)$ be a homomorphism, $A\subseteq S$, and $r\in\mathbb{N}$. If $A$ is piecewise syndetic in $S$ and $\varphi(S)$ is piecewise syndetic in $T$, then $\varphi(A)$ is piecewise syndetic in $T$.
\end{theo}

If we set $T=S$ and $\varphi(x)=xs$ for some $s\in S$, we know that $\varphi$ maintains the property of piecewise syndetic if $Ss$ is piecewise syndetic in $S$.

Next we define the group of largeness including ${\bf IP}_r^*$ that appeared in the introduction. One may see that we require some of the sets to be multisets, which is necessary in some of our proofs. Actually, if we need some of the properties of ${\bf IP}_r$, e.g. the image of ${\bf IP}_r$ sets under homomorphisms are also ${\bf IP}_r$, then mutiset is needed.

\begin{deff}
Suppose $(S,\circ)$ is a semigroup and a multiset $A$ with elements in $S$.
\begin{itemize}
\item Given a total ordering of the index set $I$ of $A=\{a_i\}_{i\in I}$, We denote ${\bf FP}(A)$ the finite product of $A$ as
\[
{\bf FP}(A)=\bigcup_{F\subseteq I\atop 1\leq|F|<\infty}\{a_{i_1}\circ\dots\circ a_{i_n}\},\quad \text{where}\quad F=\{i_1,\dots,i_n\},\quad i_1<\dots<i_n.
\]
If the operation is written as ``$+$'', as in semirings or sometimes in commutative settings, we may use ${\bf FS}$ instead which abbreviates for finite sum.
\item We say that $A$ is ${\bf IP}_r$ if $\exists F\in\mathcal{P}_f(S)$, $|F|=r$ such that ${\bf FP}(F)\subseteq A$.
\item We say that $A$ is ${\bf IP}_r^*$ if $\forall F\in\mathcal{P}_f(S)$, $|F|=r$,there is ${\bf FP}(F)\cap A\neq\emptyset$.
\end{itemize}
\end{deff}

Finally we need the concept of combinatorially rich, which is defined originally by Bergelson and Glasscock \cite{bg} in commutative semigroups and extended to arbitrary semigroups by Hindman et al \cite{hhst}. We will only need the commutative version here, and we will discuss more on it at the end of the paper.

\begin{deff}
Let $(S,+)$ be a commutative semigroup. A subset $A\subseteq S$ is combinatorially rich if for all $n\in\mathbb{N}$, there exists an $r\in\mathbb{N}$ such that for all $M\in S^{r\times n}$, there exists a non-empty $\alpha\subseteq\{1,\dots,r\}$ and $s\in S$ such that
\[
s+M_{\alpha,j}\in A,\quad \forall j\in\{1,\dots,n\}.
\]
\end{deff}

The property {\it combinatorially rich}, {\it piecewise syndetic} and {\it central} are partition regular, i.e. for any finite partition of a set with the property, one of the subsets also have the property. See, for instance, \cite{hindman}.

The concepts we have defined have deep connections between each other, and we will mention the connections that will be used in this paper. The following theorem is a property of combinatorially rich sets depicted by Bergelson and Glasscock \cite{bg}. Though the definitions on ${\bf IP_r^*}$ sets are different, the proof is the same.

\begin{theo}\label{cr}
Let $(S,+)$ be a commutative semigroup and $A\subseteq S$. Then $A$ is combinatorially rich if and only if For all $n\in N$, there exists $r\in N$ such that for all commutative semigroups $T$ and all homomorphisms $\varphi_1,\dots, \varphi_n$ from $T$ to $S$, the set
\[
\left\{d\in S:\bigcap_{i=1}^{n}(A-\varphi_i(d))\neq\emptyset\right\}
\]
is ${\bf IP}_r^*$ in $T$.
\end{theo}

It can be derived from the definition that thick sets are central, central sets are piecewise syndetic, and syndetic sets are exactly the complement of non-thick sets. We also need the following theorem, see, for instance, \cite{bg}.

\begin{theo}\label{pscr}
Let $(S,+)$ be a commutative semigroup and $A\subset S$ is piecewise syndetic. Then $A$ is combinatorially rich. 
\end{theo}

For semiring with two operations, there are some links between largeness with respect to each operation. It was shown in \cite{bh} and \cite{br} that in $\mathbb{N}$, multiplicatively central sets are additively ${\bf IP}_r$ for every $r$. Bergelson and Glasscock \cite{bg} made a generalisation to multiplicatively piecewise syndetic sets on arbitrary semirings.

\begin{theo}[Bergelson \& Glasscock]\label{m-a}
Let $(S,+,\cdot)$ be a semiring. Suppose that $(R,\cdot)$ is a subsemigroup of $(S, \cdot)$ which is additively large in the following way: $R$ is an ${\bf IP}_r$ set in $(S,+)$ for all $r$. If $A\subseteq R$ is piecewise syndetic in $(R,·)$, then $A$ is ${\bf IP}_r$ in $(S,+)$ for all $r$.
\end{theo}

In particular, setting $R=S$, we have that multiplicatively piecewise syndetic sets are additively ${\bf IP}_r$ for every $r$.

\section{Variant of the Szemer\'edi Theorem}

In this section we are to prove a piecewise syndetic form of the ${\bf IP}_r^*$ Szemer\'edi Theorem. Recall that we have defined in the introduction that
\[
D(A;\varphi_1,\dots, \varphi_n)=\left\{d\in S:\bigcap_{i=1}^{n}(\varphi_i(d)^{-1}A)\quad \text{is piecewise syndetic}\right\}.
\]
First we need the following finite Szemer\'edi's Theorem, which is from Theorem \ref{fk} together with a standard compactness argument.

\begin{lemm}\label{finite}
Let $S=(S,\cdot)$ and $R=(R,\cdot)$ be semigroups such that $R$ admits a left invariant mean. For any $n$, $k$ and homomorphisms $\varphi_1,\dots, \varphi_n$ from $S$ to $R$ such that $\varphi_i\varphi_j=\varphi_j\varphi_i$ for all $i,j$, there is $G\in\mathcal{P}_f(S)$ and $H\in\mathcal{P}_f(R)$ such that for any $k$-partition $H=\bigcup_{j=1}^kC_j$, there is $d\in G$ and $j\in [k]$ such that $\bigcap_{i=1}^{n}\left(\varphi_i(d)^{-1}C_j\right)\neq\emptyset$.
\end{lemm}
\proof We first prove the case without the finiteness limitation, i.e. the case when $G=S$ and $H=R$. Pick an $s\in S$ arbitrarily. Denote transformations
\[
T_i(a)=\varphi_i(s)a,\quad 1\leq i\leq n.
\]
Then
\[
T_i^{-1}(A)=\varphi_i(s)^{-1}A
\]
and thus $T_i$'s are measure preserving, and
\[
T_iT_j(a)=\varphi_j(s)\varphi_i(s)a=\varphi_i(s)\varphi_j(s)a=T_jT_i(a),
\]
i.e. $(T_i)_{i=1}^n$ are commuting. By Theorem \ref{fk}, there is some $j\in[k]$ and $m\in\mathbb{N}$ such that
\[
\bigcap_{i=1}^n\left(\varphi_i\left(s^m\right)^{-1}C_j\right)=\bigcap_{i=1}^n T_n^{-m}C_j\neq\emptyset.
\]

Therefore, for any $k$-partition $R=\bigcup_{j=1}^kC_j$, we can find $d=s^m\in S$ and $j\in [k]$ such that $\bigcap_{i=1}^{n}\left(\varphi_i(d)^{-1}C_j\right)\neq\emptyset$. We identify each $k$-partition to a function (element) in $[k]^R$, equipped with the product topology ($[k]=\{1,\dots,k\}$ equipped with the discrete topology). By the Tychonoff Theorem, $[k]^R$ is compact. For $d\in S$, $a\in R$, Define
\[
\mathcal{C}(d,a)=\left\{c\in[k]^R:\forall j\in[k],a\notin\bigcap_{i=1}^{n}\left(\varphi_i(d)^{-1}c^{-1}(j)\right)\right\},
\]
(here $c^{-1}(j)$ is identical with $C_j$.) One can see it is closed in $[k]^R$. What we have already proved says that
\[
\bigcap_{(d,a)\in S\times R}\mathcal{C}(h)=\emptyset,
\]
and by the property of compactness, there is some $I\in\mathcal{P}_f(S\times R)$ such that
\[
\bigcap_{(d,a)\in I}\mathcal{C}(d,a)=\emptyset.
\]
Find $G\in\mathcal{P}_f(S), H'\in\mathcal{P}_f(R)$ such that $I\subseteq G\times H'$. Then for all $k$-partition $c\in[k]^R$, there is a $d\in G$ and $j\in [k]$ such that
\[
H'\cap\bigcap_{i=1}^{n}\left(\varphi_i(d)^{-1}c^{-1}(j)\right)\neq\emptyset.
\]
Set
\[
H=\bigcup_{d\in G}\bigcup_{i\in[n]}\varphi_i(d)H',
\]
then
\[
H'\subseteq\bigcap_{i=1}^{n}\left(\varphi_i(d)^{-1}H\right),
\]
hence
\[
\bigcap_{i=1}^{n}\left(\varphi_i(d)^{-1}H\right)\cap\bigcap_{i=1}^{n}\left(\varphi_i(d)^{-1}c^{-1}(j)\right)=\bigcap_{i=1}^{n}\left(\varphi_i(d)^{-1}\left(H\cap c^{-1}(j)\right)\right)\neq\emptyset.
\]
This completes the proof.\qed

Now we prove the first half of the piecewise syndetic Szemer\'edi Theorem, which is a generalisation of the lemma above, and is also needed while proving the second half. The basic idea is that if there is $s$ such that the intersection is not empty, then with the property of partition regularity we can require the intersection to be piecewise syndetic.

\begin{theo}\label{noncommu}
Let $S=(S,\cdot)$ and $R=(R,\cdot)$ be semigroups such that $R$ admits a left invariant mean, and $A\subseteq R$ be piecewise syndetic. Then for any $n$ and commuting homomorphisms $\varphi_1,\dots, \varphi_n$ from $S$ to $R$, there is some $t\in R$ such that
\[
D(t^{-1}A;\varphi_1,\dots, \varphi_n)\neq\emptyset.
\]
\end{theo}
\proof We add, if there is none, an identical element $\bm{1}$ to the semigroup $R$, that is, $\bm{1}r=r\bm{1}$ for all $r\in R$, and denote $R'=R\cup\{\bm{1}\}$. Though $R$ may have identical element (left or right) itself, there is no contradiction because it need not to be identical in $R'$.

We claim that $A$ is still piecewise syndetic in $R'$. Recall that a set is piecewise syndetic if it is the intersection of a syndetic set and a thick set, and a set is thick if and only if its complement is syndetic. Thus it suffices to prove that if $M$ is syndetic in $R$, then $M$ (and hence $M\cup\{\bm{1}\}$ as well) is syndetic in $R'$. But that is trivial, as we can write by definition
\[
R=\bigcap_{i=1}^ms_m^{-1}M,
\]
and pick $a\in M$, then $\bm{1}\in a^{-1}M$, hence
\[
R'=a^{-1}M\cap\bigcap_{i=1}^ms_m^{-1}M.
\]

We also add a homomorphism $\varphi_0(s)=\bm{1}$ for all $s\in S$ if there is none. Then it suffices to prove that
\[
D(A;\varphi_0,\varphi_1,\dots, \varphi_n)\neq\emptyset.
\]

Therefore, we may suppose without loss of generality that, $R$ bears a identical element $\bm{1}$ and that $\varphi_1(s)=\bm{1}$ for all $s\in S$. 

From the definition, we know that there is $F\in\mathcal{P}_f(R)$ such that
\[
T=\bigcup_{t\in F} t^{-1}A
\]
is thick in $R$. Find $G$ and $H$ that satisfy the conditions in Lemma \ref{finite}, then there is some $x\in R$ such that
\[
Hx\subseteq T.
\]

Note that when considering the set $\bigcap_{i=1}^{n}\left(\varphi_i(s)^{-1}C_j\right)$, it does not affect if we right multiply the sets with a constant $x$, and hence if $H$ satisfies the condition in Lemma \ref{finite}, then so is $Hx$ as well. Since
\[
Hx=\bigcup_{t\in F}\left(Hx\cap\left(t^{-1}A\right)\right),
\]
there is some $d\in G$ and $t\in F$ such that
\[
\bigcap_{i=1}^{n}\left(\varphi_i(d)^{-1}\left(Hx\cap\left(t^{-1}A\right)\right)\right)\neq\emptyset.
\]
Next we consider the set
\[
A'=\left\{a\in R:\exists t\in F,\exists d\in G, a\in\bigcap_{i=1}^{n}\left(\varphi_i(d)^{-1}t^{-1}A\right)\right\}.
\]
We claim that
\[
\bigcup_{h\in H'}h^{-1}A'
\]
is thick, which implies that $A'$ is piecewise syndetic. Following the definition, for any $F'\in\mathcal{P}_f(R)$, the set $HF'$ is also finite, so by the thickness of $T$, there is some $x$ such that
\[
HF'x\subseteq T.
\]
We are going to show that
\[
F'x\in\bigcup_{h\in H'}h^{-1}A',
\]
which suffices if for all $c\in F'$, there is some $a\in A'$ and $h\in H'$ such that $cx=h^{-1}a$. Note that $G$ and $Hc$ also satisfy the conditions in Lemma \ref{finite}, so by the same analysis, there is some $d\in G$ and $t\in F$ such that
\[
\bigcap_{i=1}^{n}\left(\varphi_i(d)^{-1}\left(Hcx\cap t^{-1}A\right)\right)\neq\emptyset.
\]
Pick $a$ from the above intersection. Then $a\in A'$. Moreover, since $\varphi_1(d)=e$ by the assumption in the beginning, we also know that $a\in Hcx$, which implies that there is some $h\in H'$ such that $cx=h^{-1}a$. Therefore, $A'$ is piecewise syndetic.

Finally, since
\[
A'=\bigcup_{t\in F}\bigcup_{d\in G}\bigcap_{i=1}^{n}\left(\varphi_i(d)^{-1}A\right),
\]
from the finiteness of $F$,$G$ and the partition regularity of piecewise syndetic property, there is some $t\in F$ and $d\in G$ such that $\bigcap_{i=1}^{n}\left(\varphi_i(d)^{-1}t^{-1}A\right)$ is piecewise syndetic, i.e.
\[
D(t^{-1}A;\varphi_1,\dots, \varphi_n)\neq\emptyset.
\]
This completes the proof.\qed

Note that one cannot require $D(A;\varphi_1,\dots, \varphi_n)$ to be nonempty. For example, let $R$ be the free semigroup generated by two elements $x$ and $y$, $S$ be the sub-semigroup generated by only $x$ (without the empty word), and $\phi$ be the identical homomorphism. We set $A=yR$, which is syndetic in $R$, but for all $s\in S$, the set $\phi(s)^{-1}A$ is empty.

Now we are ready to prove the second half, which is based on the proof of the above Theorem together with properties of combinatorially rich sets. We will omit some details.

\begin{theo}\label{commu}
Let $S=(S,+)$ and $R=(R,+)$ be commutative semigroups and $A\subseteq R$ be piecewise syndetic. Then for any $n$, there exists $r$ such that
\[
D(A;\varphi_1,\dots, \varphi_n)
\]
is ${\bf IP}_r^*$ in $S$ for all homomorphisms $\varphi_1,\dots, \varphi_n$ from $S$ to $R$.
\end{theo}
\proof By Theorem \ref{pscr}, $A$ is combinatorially rich in $R$. Then Theorem \ref{cr} states that, the set of $d$ such that
\[
\bigcap_{i=1}^n(A-\varphi_i(d))\neq\emptyset
\]
is ${\bf IP}_r^*$ in $S$, which implies that for every ${\bf IP}_r$ set $G\subseteq S$ and any finite partition $R=\bigcup_{j=1}^kC_j$, some $C_j$ is combinatorially rich in $R$ by partition regularity, and in turn there is $d\in G$ such that
\[
\bigcap_{i=1}^{n}\left(\varphi_i(d)^{-1}C_j\right)\neq\emptyset.
\]
By a similar compact argument as in the proof of Lemma \ref{finite}, for every finite ${\bf IP}_r$ set $G\subseteq S$, there is $H\in\mathcal{P}_f(R)$ such that for any finite partition of $H$, the above is also true. Following the proof of Theorem \ref{noncommu}, for every such $G$, there is some $t\in F$ and $d\in G$ such that the intersection of all $A-t-\varphi_i(d)$ is not only nonempty but also piecewise syndetic in $R$. From the definition, we know that
\[
\bigcap_{i=1}^n(A-\varphi_i(d))=\bigcap_{i=1}^n(A-t-\varphi_i(d))+t
\]
is piecewise syndetic in $R$ for such $d\in G$, and hence
\[
D(A;\varphi_1,\dots, \varphi_n)\cap G\neq\emptyset.
\]
From the definition, we know that $D(A;\varphi_1,\dots, \varphi_n)$ is ${\bf IP}_r^*$ in $S$.\qed

\begin{remark}\label{remark}
In the definition of $D$, we used the expression $\bigcap_{i=1}^{n}(\varphi_i(d)^{-1}A)$. However, as we stated in the proof of Theorem \ref{noncommu}, since adding an identical element does not affect the property of piecewise syndetic, one can change the expression into
\[
A\cap\bigcap_{i=1}^{n}(\varphi_i(d)^{-1}A),
\]
which is sometimes the form we need in practise.
\end{remark}

\section{Products and sums in semirings}

To extend Bowen's results to semirings, we need to make use of following lemma that Bowen proved in \cite{bowen}. Although Bowen proved it only on $\mathbb{N}$, the proof do not need the commutative law for multiplication and hence holds on arbitrary semiring. We need to modify the expression a bit to fit for the semiring setting. As naturally defined, for two subsets $A$ and $B$ in semigroup $(S,\cdot)$, their product $AB$ represents the set $\{ab :a\in A,b\in B\}$.

\begin{theo}[Bowen]\label{keylemma}
Let $S=(S,+,\cdot)$ be a semiring, $L\subseteq (\beta S,\cdot)$ be a minimal left ideal, $e\in L$ be an idempotent,
and $p\in L$. Then for any $A\in p$, $B\in e$, and $f:\omega\to\omega$ there is a sequence of sets $F_0,F_1,...\subseteq S$ with $|F_i|=f(i)$ such that for any $n$ and $I=\{i_1,\dots,i_n\}$, the set
\[
{\bf FS}(F_{i_0}){\bf FS}(F_{i_1})\dots{\bf FS}(F_{i_n})
\]
is in $A$ if $0\in I$ and is in $B$ otherwise.
\end{theo}

We shall illustrate what it mean if we just want two sets $F_0,F_1$, which is the case we need if we want to find a monochromatic $\{x,y,x+y,xy\}$. By standard result of ultrafilter method, one can find as large as desired (but finite) subsets $F_0,F_1$ such that ${\bf FS}(F_0)\subseteq A$ and ${\bf FS}(F_1)\subseteq B$ (recall that multiplicatively piecewise syndetic means additively ${\bf IP}_r$ for all $r$, see, for instance, \cite{bg}). This theorem tells us that, one can in addition require that the set ${\bf FS}(F_0){\bf FS}(F_1)$ is also contained in $A$, which is the key ingredient in finding $\{x,y,x+y,xy\}$ patterns.

We split the result into two theorems, whose proofs are both based on Bowen's proofs in \cite{bowen}. 

\begin{theo}
Let $S=(S,+,\cdot)$ be a semiring such that $Ss$ is piecewise syndetic in $(S,+)$ for all $s\in S$. Let $T$ be thick in $(S,\cdot)$, $M_0,M_1$ be syndetic in $(S,\cdot)$, and $S=C_0\cup C_1$ be a bipartition such that $C_0\in T\cap M_0$ and $C_1\in T\cap M_1$. Then for all $k\in\mathbb{N}$, there is some $y\in S$ and $i\in\{0,1\}$ such that
\[
\Big\{x:\{x,y,xy,x+y,x+2y,\dots,x+ky\}\subseteq C_i\Big\}
\]
is also piecewise syndetic in $(S,+)$.
\end{theo}
\proof The assupmtion ensures that there is minimal left ideal $L\subseteq\beta S$ and idempotent $e\in L$ and $p\in L$ such that $C_0\in p\in L$ and $C_1\in e\in L$. Then the condition of Theorem \ref{keylemma} is satisfied and there is as large as desired sets $F_0,F_1\in\mathcal{P}_f(S)$ such that
\[
{\bf FS}(F_0)\subseteq C_0,\quad {\bf FS}(F_1)\subseteq C_1,\quad {\bf FS}(F_0){\bf FS}(F_1)\subseteq C_0.
\]
We first consider the case when $C_0$ is not additively piecewise syndetic. Then by the partition regularity, we know that $C_1$ is piecewise syndetic in $(S,+)$. Using Theorem \ref{commu} and Remark \ref{remark} on $C_1$ and $F_1$, we find $y\in {\bf FS}(F_1)$ such that
\[
A=\bigcap_{i=0}^{k}(C_1-iy)
\]
is piecewise syndetic in $(S,+)$. By Theorem \ref{homo}, $Ay$ is piecewise syndetic in $(S,+)$. Since $C_0$, and hence $C_0\cap Ay$ is not piecewise syndetic, we know that $Ay\cap C_1$ is piecewise syndetic. For $y$ and all $x'\in Ay\cap C_1$,
\[
\left\{x'y^{-1},y,x',x'y^{-1}+iy:i\in[k]\right\}\in C_1.
\]
Moreover, by Theorem \ref{homo} we know that
\[
x=x'y^{-1}\in(Ay\cap C_1)y^{-1},
\]
which is also piecewise syndetic in $(S,+)$, and hence we are done.

Otherwise, $C_0$ is additively piecewise syndetic, and we find $y_0\in {\bf FS}(F_0)$ such that
\[
A_0=\bigcap_{i=0}^{k}(C_0-iy_0)
\]
is piecewise syndetic in $(S,+)$. If $A_0y_0\cap C_0$ is piecewise syndetic then we are done as before. Otherwise, $C_1’=A_0y_0\cap C_1$ is piecewise syndetic in $(S,+)$. Again using Theorem \ref{commu} and Remark \ref{remark} on $C_1'$ and $F_1$, we find $y_1\in {\bf FS}(F_1)$ such that
\[
A_1=\bigcap_{i=0}^{k}(C_1'-iy_1)\cap\bigcap_{i=1}^{k}(C_1'-iy_0y_1y_0)
\]
is piecewise syndetic in $(S,+)$, and hence $A_1y_1$ is piecewise syndetic in $(S,+)$. If $A_1y_1\cap C_1$ is piecewise syndetic then again we are done as before. Otherwise, $A_1y_1\cap C_0$ is piecewise syndetic in $(S,+)$. Here,
\[
y=y_0y_1\in{\bf FS}(F_0){\bf FS}(F_1)\subseteq C_0.
\]
For all
\[
x'\in A_1y_1\cap C_0\subseteq C_0,
\]
we have
\[
x=x'y^{-1}\in (A_1y_1\cap C_0)y^{-1}\subseteq C_1'y_0y_1y^{-1}=A_0'\subseteq C_0,
\]
and for all $0\leq i\leq k$,
\[
x+iy=x'y_1^{-1}y_0^{-1}+iy_0y_1y_0y_0^{-1}\in A_1y_0^{-1}\subseteq C_0.
\]
Similarly, the set of $x$ that satisfies the condition is piecewise syndetic.\qed

We have proved the case when $C_0$ and $C_1$ are both multiplicatively syndetic, whose reversed side is that one of the $C_i$'s is not syndetic, i.e. the other is multiplicatively thick.

\begin{theo}
Let $S=(S,+,\cdot)$ be a semiring, $S=C_0\cup C_1$ be a bipartition such that $C_0$ is multiplicatively thick. Then for all $k,l\in\mathbb{N}$, there is $i\in\{0,1\}$ and $x,y\in S$ such that
\[
\{x,y,kx+y,xy,x^2y,\dots,x^ly\}\subseteq C_i.
\]
\end{theo}
\proof For every $n$, we denote $F[x_1,\dots,x_n]$ to be the free semigroup on $n$ elements, i.e. elements in $F[x_1,\dots,x_n]$ are words of finite length $(a_1,a_2,\dots,a_m)$ such that $a_j\in\{x_1,\dots,x_n\}$ for all $1\leq j\leq m$. We denote $F_2[x_1,\dots,x_n]$ to be the subset of $F$ with all words such that $x_i$ appears no more than $2l^2$ times for all $1\leq i\leq n$, and $F_{2,k}[x_1,\dots,x_n]$ all linear combinations with elements in $F_2[x_1,\dots,x_n]$ and coefficients in $\mathbb{N}$ which sum up to no more than $k^{3nkl^2}$.

The set $F_{2,k}[x_1,\dots,x_n]$ is finite, and since $C_0$ is multiplicatively thick, we can find $a_1,a_2,a_3,a_4,a_5$ such that $a_1\in C_0$, and
\[
F_{2,k}[a_1,\dots,a_{i-1}]a_i\in C_0,\quad \forall\ 2\leq i\leq5.
\]

Suppose by contrary that there are no such patterns. Since
\[
\{ka_1,ka_2,k^2a_1a_2,\dots,k^{l+1}a_1^la_2\}\subseteq C_0,
\]
there must be
\[
k^2a_1+ka_2\in C_1.
\]
Similarly, $ka_3+a_4,ka_4+a_5\in C_1$. Moreover, for $x=\left(k^2a_1+ka_2\right)^ia_3$ and $y=\left(k^2a_1+ka_2\right)^ia_4$, $1\leq i\leq l$, we have $x,y\in A$, and
\[
x^jy=\left(k^2a_1+ka_2\right)^{i(j+1)}a_3^ja_4\in C_0,\quad \forall\ 1\leq j\leq l,
\]
hence
\[
kx+y=\left(k^2a_1+ka_2\right)^i(ka_3+a_4)\in C_1,\quad \forall\ 1\leq i\leq l.
\]
For $x=k^2a_1+ka_2$ and $y=ka_3+a_4$, we have illustrated that $x,y\in C_1$ and for all $1\leq i\leq l$, $x^iy\in C_1$, hence
\[
kx+y=k^3a_1+k^2a_2+ka_3+a_4\in C_0.
\]
For $x=k^2a_1+ka_2+a_3$ and $y=a_4\in C_0$, we have $kx+y\in C_0$ by above, and $x^iy\in C_0$ for all $1\leq i\leq l$ by the way we find $a_4$. Hence
\[
x=k^2a_1+ka_2+a_3\in C_1.
\]
For $x=k^3a_1+k^2a_2+ka_3+a_4$ and $y=a_5$, they are both in $C_0$ and $x^iy\in C_0$ for all $1\leq i\leq l$. Hence
\[
kx+y=k^4a_1+k^3a_2+k^2a_3+ka_4+a_5\in C_1.
\]
For $x=k^2a_1+ka_2+a_3$ and $y=ka_4+a_5$, we have showed that both of them are in $C_1$, and $kx+y\in C_1$. Hence there is some $i\in[l]$ such that
\[
x^iy=\left(k^2a_1+ka_2+a_3\right)^i(ka_4+a_5)\in C_0.
\]
Now for $x=\left(k^2a_1+ka_2+a_3\right)^ia_4$ and $y=\left(k^2a_1+ka_2+a_3\right)^ia_5$, we have
\[
\{x,y,kx+y,xy,x^2y,\dots,x^ly\}\subseteq C_i,
\]
a contradiction.\qed

The above are Bowen's results on 2 variants. The result on $n$ variants can also be extended to semirings. The proof is similar to the original proof in \cite{bowen} together with our discussion to fit for semiring settings.

\begin{theo}
Given any semiring $(S,+,\cdot)$  such that $Ss$ is piecewise syndetic in $(S,+)$ for all $s\in S$, any $2$-coloring of $S$ and $n\in\mathbb{N}$, there are $x_1, ..., x_n\in S$ such that
\[
\{x_i,x_1x_2\dots x_i,\sum_{j=1}^nx_j:i\leq n\}
\]
is monochromatic.
\end{theo}

The above theorems all rely on discussions on both color class, and we cannot expect on finding the wanted pattern in a certain subset, which seems to be the core difficulty when considering any finite partition rather than only 2-partitions. If we use the thick property in a different way, we may obtain the following theorem.

\begin{theo}
For any $n\in\mathbb{N}$, any semiring $S=(S,+,\cdot)$ and $A\subseteq S$, if $A$ is both multiplicatively thick and syndetic, then there is $x_1,\dots,x_n$ such that
\[
\{x_i,x_1x_2\dots x_n,\sum_{i=1}^{n}a_ix_i:1\leq i\leq n, 0\leq a_i\leq k\}\subseteq A
\]
\end{theo}
\proof Since $A$ is syndetic in $(S,\cdot)$, we can find $s_1,s_2,\dots,s_m\in S$ such that
\[
s_1^{-1}A\cap s_2^{-1}A\cap\dots\cap s_m^{-1}A=S.
\]
For a set $A$, denote $^kA$ the multiset that each element in $A$ appears $k$ times in $^kA$. Pick arbitrarily $s_{m+1},\dots,s_{m+n-1}\in S$, and denote
\[
F={\bf FS}(^k\{s_1,\dots,s_m,\dots,s_{m+n-1}\}),
\]
Then $F\in\mathcal{P}_f(S)$. Since $A$ is thick in $(S,\cdot)$, we can find $x\in S$ such that
\[
Fx\subseteq A
\]
Now we consider the element $xs_{m+1}xs_{m+2}x\dots s_{m+n-1}x$. Since $A$ is syndetic, we know that there is some $i\in[n]$ such that
\[
s_ixs_{m+1}xs_{m+2}x\dots s_{m+n-1}x\in A.
\]
Set $x_1=s_ix$, $x_j=s_{m+j-1}x$ for $2\leq j\leq n$, then $x_i$ are in $A$, there product is in $A$, and there linear combination is in $A$ from $Fx\subseteq A$ and the definition of $F$.\qed

\section{Concluding remarks}

We first proved in section 3 that given suitable semigroups $(S,\cdot)$ and $(R,\cdot)$ and homomorphisms $\varphi_1,\dots,\varphi_n$, for all piecewise syndetic $A\subseteq R$ one can find $d\in S$ and a piecewise syndetic set $a\in R$ such that $\{a+\varphi_i(d):1\leq i\leq n\}$ is in $t^{-1}A$ for some $t\in R$. We also showed that if the semigroups are commutative, then one can find a ${\bf IP}_r^*$ set of $d$ that satisfies the property. We used the property of combinatorially rich which was originally defined in commutative groups.

However, Hindman et al \cite{hhst} has extended the concept of combinatorially rich to arbitrary semigroups. The definition agrees on the commutative case, and combinatorially rich sets in arbitrary semigroups have similar properties with those in commutative semigroups. For instance, piecewise syndetic sets are combinatorially rich. However, the method in proving Theorem \ref{cr} is not valid in the general setting, and we wonder whether there are other methods to deal with the problem.

\begin{conj}
Let $S=(S,\cdot)$ and $R=(R,\cdot)$ be semigroups such that $R$ admits a left invariant mean, and $A\subseteq R$ be piecewise syndetic. Then for all $n$, there is $t\in R$ and $r\in\mathbb{N}$ such that
\[
D(t^{-1}A;\varphi_1,\dots, \varphi_n)
\]
is ${\bf IP}_r^*$ in $S$ for all commuting homomorphisms $\varphi_1,\dots, \varphi_n$ from $S$ to $R$.
\end{conj}

Especially, if Theorem \ref{cr} still holds for arbitrary semigroup, then the conjecture is true with the same method in proving Theorem \ref{commu}.

In section 4, we extended Bowen's result to semirings $(S,+,\cdot)$ that $Ss$ is piecewise syndetic for all $s\in S$, showing that essentially, one can find a monochromatic $\{x,y,x+y,xy\}$ in all bipartition of $S$. If $S$ does not satisfy the property itself but some subsemiring $R\subseteq S$ satisfies, then the method is also valid. It is natural to ask what the case is when there is no such $R$. An example is the polynomial semiring $\mathbb{N}[x]$.

\begin{conj}
For all bipatition of the polynomial semiring $\mathbb{N}[x]=C_0\cup C_1$, there is some $i\in\{0,1\}$ and $f,g\in\mathbb{N}[x]$ such that $\{f,g,f+g,fg\}\subseteq C_i$.
\end{conj}

\noindent

\noindent{\large\bf Acknowledgement.}  

Thanks to Hehui Wu for the early discussion of this project
and the helpful comments.

\end{document}